\documentclass[10pt]{amsart}

\usepackage{amsthm,amsmath,amssymb,amscd,amsfonts,latexsym}

\theoremstyle{plain}
\newtheorem{theorem}{Theorem}[section]

\newtheorem{proposition}[theorem]{Proposition}
\newtheorem{maintheorem}[theorem]{Main Theorem}
\newtheorem{criterion}[theorem]{Criterion}

\newtheorem{def-thm}[theorem]{Definition-Theorem}
\newtheorem{lemma}[theorem]{Lemma}

\theoremstyle{definition}
\newtheorem{definition}[theorem]{Definition}

\newtheorem{remark}[theorem]{Remark}

\newtheorem{facts}[theorem]{Facts}

\newtheorem*{acknowledgement}{Acknowledgement}

\newcommand{\PP}{\mathbb{P}}

\newcommand{\NN}{\mathbb{N}}
\newcommand{\ZZ}{\mathbb{Z}}
\newcommand{\CC}{\mathbb{C}}

\newcommand{\OO}{{\mathcal O}}

\newcommand{\II}{{\mathcal I}}

\newcommand{\mm}{{\mathfrak m}}

\newcommand{\fra}{{\mathfrak a}}

\newcommand{\ke}{{K\"ahler-Einstein\ }}
\newcommand{\kr}{{K\"ahler-Ricci\ }}

\newcommand{\ddt}{\frac{\partial}{\partial t}}

\DeclareMathOperator{\Dom}{Dom}
\DeclareMathOperator{\Pic}{Pic}

\DeclareMathOperator{\Aut}{Aut}
\DeclareMathOperator{\Cr}{Cr}
\DeclareMathOperator{\id}{id}
\DeclareMathOperator{\rank}{rank}

\DeclareMathOperator{\PGL}{PGL}

\begin{document}

\title[Convergence of the K\"ahler-Ricci flow and multiplier ideal sheaves]
{Convergence of the K\"ahler-Ricci flow and multiplier ideal sheaves on del Pezzo surfaces}
\begin{abstract} 
On certain del Pezzo surfaces with large automorphism groups, it is shown that the solution to the K\"ahler-Ricci flow with a certain initial value converges in $C^\infty$-norm exponentially fast to a K\"ahler-Einstein metric. The proof is based on the method of multiplier ideal sheaves.
\end{abstract}

\author[G. Heier]{Gordon Heier}
\address{Department of Mathematics\\ Christmas-Saucon Hall\\ Lehigh University\\ \ \ 14 E.~Packer Avenue \\
Bethlehem, PA 18015\\ USA}

\email{heier@lehigh.edu}


\maketitle

\section{Introduction} 
Let $X$ be an $n$-dimensional compact complex manifold with positive first Chern class $c_1(X)$. Such manifolds are called {\it Fano} manifolds. The \kr flow on $X$ is defined by the equation
\begin{equation}\label{RF}
\ddt g_{i\bar{j}}=-R_{i\bar{j}}+g_{i\bar{j}},
\end{equation}
where $R_{i\bar{j}}=-\partial_i\partial_{\bar{j}} \log\det g_{\alpha\bar\beta}$ is the Ricci curvature tensor of the hermitian metric $\sum_{i,j} g_{i\bar{j}}dz_i\otimes d\bar{z}_{j}$. If the class of the K\"ahler form $\hat{\omega}=\frac{i}{2\pi}\sum_{i,j} \hat{g}_{i\bar{j}}dz_i\wedge d\bar{z}_{j}$ is $c_1(X)$, then the \kr flow preserves the class of $i\sum_{i,j} \hat{g}_{i\bar{j}}dz_i\wedge d\bar{z}_{j}$, so we can write
\begin{equation*}
g_{i\bar{j}}=\hat{g}_{i\bar{j}}+\partial_i\partial_{\bar{j}}\phi,
\end{equation*}
for the solution to the \kr flow with initial condition
\begin{equation*}
g_{i\bar{j}}(0)=\hat{g}_{i\bar{j}}.
\end{equation*}\par
Equation \eqref{RF} can be reformulated to
\begin{equation}\label{RF_phi}
\ddt \phi=\log\frac{\det g_{\alpha\bar{\beta}}}{\det \hat{g}_{\alpha\bar{\beta}}}+\phi-\hat f,\quad \phi(0)=c_0\in \CC,
\end{equation}
where $\hat f$ is the Ricci potential, ie for $\hat {R}_{i\bar{j}}=-\partial_i\partial_{\bar{j}} \log\det \hat{g}_{\alpha\bar\beta}$, we have $\hat{R}_{i\bar{j}}-\hat{g}_{i\bar{j}}=\partial_i\partial_{\bar{j}} \hat f$.
It was proven in \cite{Cao_Inv_Math_1985} that the solution to \eqref{RF} exists for all $t>0$. The present paper investigates the issue of convergence, based on the following Theorem which first appeared \cite{PSS}. The version given below, which is stronger than the one in \cite{PSS}, is based on \cite{PS_CDM}.

\begin{theorem}[\cite{PSS,PS_CDM}]\label{PSS_Thm}
Let $X$ be a Fano manifold. Consider the Ricci flow in the form of \eqref{RF_phi} with the initial value $c_0$ specified by \cite[(2.10)]{PSS}. The following two statements are equivalent.
\begin{enumerate}
\item There exists $p>1$ such that
\begin{equation*}
\sup_{t\geq 0} \int_X e^{-p\phi} \hat{\omega}^n<\infty.
\end{equation*}
\item The family of metrics $g_{i\bar{j}}(t)$ converges in $C^\infty$-norm exponentially fast to a \ke metric.
\end{enumerate}
\end{theorem}
The preceding theorem will allow us to formulate the sufficient Criterion \ref{simple_crit} for the above statement (ii) to hold, in analogy to Nadel's criterion for the existence of \ke metrics (see \cite{Nadel_Annals, Demailly_Kollar, Heier_ke_del_pezzo}). It is well-known that some Fano manifolds do not possess a \ke metric (eg $\PP^2$ blown up in one or two points, see Section \ref{low_degree}), so we cannot expect (ii) to hold in general on a Fano manifold. In this paper, we mention no necessary condition for (ii) other than the existence of a \ke metric.

\par
First, we quickly recall the basics of multiplier ideal sheaves. The following is the standard definition of the multiplier ideal sheaf pertaining to a plurisubharmonic function on a complex manifold.
\begin{def-thm}[\cite{Nadel_Annals}]
Let $\varphi$ be a plurisubharmonic function on the complex manifold $X$. Then the {\it multiplier ideal sheaf} $\II(\varphi)$ is the subsheaf of $\OO_X$ defined by
\begin{equation*}
\II(\varphi)(U) =\{f\in \OO_X(U):|f|^2e^{-\varphi}\in L^1_{\text{loc}}(U)\}
\end{equation*}
for every open set $U\subseteq X$. It is a coherent subsheaf.
\end{def-thm}
Multiplier ideal sheaves have turned out to be very useful in algebraic geometry, mainly because of the following vanishing theorem. They are usually defined using the notion of a singular hermitian metric on a line bundle, which in general is a metric $h$ that is given on a small open set $U$ by $h=e^{-\varphi}$, where $\varphi$ is $L^1(U)$. If $\varphi$ is plurisubharmonic for every $U$, the multiplier ideal sheaf $\II(h)$ attached to $h$ is defined by $\II(h)(U)=\II(\varphi)(U)$ if $h=e^{-\varphi}$ on $U$.
\begin{theorem}[Nadel's vanishing theorem]
Let $X$ be a compact complex K\"ahler manifold. Let $L$ be a line bundle on $X$ equipped with a singular hermitian metric such that the curvature current $-\frac{i}{2\pi}\partial\bar\partial \log h$ is positive definite in the sense of currents, ie there is a smooth positive definite $(1,1)$-form $\omega$ and $\varepsilon >0$ such that $-\frac{i}{2\pi}\partial\bar\partial \log h \geq \varepsilon\omega$. Then
\begin{equation*}
H^q(X,(K_X+L)\otimes \II(h))=0 \quad \text{for all } q\geq 1.
\end{equation*}
\end{theorem}
We now develop a Nadel-type criterion for Theorem \ref{PSS_Thm} (ii) to hold true.
If we assume that Theorem \ref{PSS_Thm} (ii) does not hold true, then, according to the theorem, for all $p>1$ there exists a sequence of times $t_i\to \infty$ with
\begin{equation*}
\lim_{i\to\infty} \int_X e^{-p\phi(t_i)} \hat{\omega}^n =\infty.
\end{equation*}
In fact, also
\begin{equation*}
\lim_{i\to\infty} \int_X e^{-p(\phi(t_i)-\frac 1 V \int_X \phi(t_i)\hat{\omega}^n)} \hat{\omega}^n =\infty,
\end{equation*}
where $V=\int_X \hat{\omega}^n$.
Let $\psi$ be an $L^1$ limit of the sequence $\phi(t_i)-\frac 1 V \int_X \phi(t_i)\hat{\omega}^n$. By semi-continuity, $||e^{-\psi}||_{L^p(X)}=\infty$. If $G\subseteq \Aut(X)$ is a compact subgroup and $\hat \omega$ is $G$-invariant, then we can assume $\psi$ and $\II(-p\psi)$ to be $G$-invariant as well.\par
We have
\begin{equation*}
\hat{\omega}+\frac{i}{2\pi}\partial\overline\partial \psi \geq 0.
\end{equation*}
Let $\hat{h}$ be a smooth $G$-invariant hermitian metric for the anticanonical line bundle $-K_X$ with $\frac{1}{2\pi i}\partial\overline\partial \log \hat{h}=\hat{\omega}\in c_1(X)$. The singular $G$-invariant hermitian metric $\hat{h}^{1+\lfloor p \rfloor}\cdot e^{-p\psi}$ is a singular metric for $-(1+\lfloor p \rfloor)K_X$ with positive curvature in the sense of currents:
\begin{eqnarray*}
-\frac{i}{2\pi}\partial\bar\partial \log (\hat{h}^{1+\lfloor p \rfloor}\cdot e^{-p\psi})&=&-\frac{i}{2\pi}\partial\bar\partial \log (\hat{h}^{1+\lfloor p \rfloor-p}\cdot \hat{h}^{p}\cdot e^{-p\psi})\\
&\geq&(1+\lfloor p \rfloor-p)\hat{\omega}.
\end{eqnarray*}
Note that for all $p>1$, we have $1+\lfloor p \rfloor-p>0$.\par
Letting $p=\frac{3}{2}$ (or any other number in the interval $]1,2[$) yields the following.
\begin{theorem}[Nadel-type criterion]\label{Nadel_type_crit}
Let $X$ be a Fano manifold. Assume that statement (ii) of Theorem \ref{PSS_Thm} does not hold. Then the $G$-invariant singular hermitian metric $h=\hat{h}^2\cdot e^{-\frac{3}{2}\psi}$ on the line bundle $-2K_X$ is such that
\begin{enumerate}
\item the curvature of $h$ is positive definite in the sense of currents,
\item $0\not =\II(\frac{3}{2}\psi)\not=\OO_X$.
\end{enumerate}
The multiplier ideal sheaf $\II(\frac{3}{2}\psi)$ is also $G$-invariant. In particular, every element of $G$ maps the zero-set $V(\II(\frac{3}{2}\psi))$ to itself.
\end{theorem}
Note that we can apply Nadel's vanishing theorem with $\tilde h=hh_E=\hat{h}^2\cdot e^{-\frac{3}{2}\psi}h_E$ and $L=-2 K_X+E$, where $E$ is an arbitrary line bundle with semi-positive metric $h_E$, to obtain
\begin{equation}\label{van_RF}
H^q(X,(K_X+L)\otimes \II(h))=H^q(X,(-K_X+E)\otimes \II(\tfrac{3}{2}\psi))=0 \quad \text{for all } q\geq 1.
\end{equation}
\begin{definition} An ideal subsheaf $\II\subseteq\OO_X$ is said to satisfy {\it Property (Van)} if for every semi-positive line bundle $E$
\begin{equation*}
H^q(X,(-K_X+E)\otimes \II)=0 \quad \text{for all } q\geq 1.
\end{equation*}
\end{definition}
The above discussion can be summed up in the following sufficient criterion.
\begin{criterion}\label{simple_crit}
Let $X$ be a Fano manifold. Let $G$ be a compact subgroup of $\Aut(X)$. Let there be no nontrivial $G$-invariant subsheaf $\II\subseteq \OO_X$ which satisfies Property (Van). Then statement (ii) of Theorem \ref{PSS_Thm} holds true.
\end{criterion}
A criterion of this kind is the essence of Nadel's technique, which can be applied under similar circumstances based on the continuity method for the Monge-Amp\`ere equation to show the existence of \ke metrics on certain Fano manifolds \cite{Nadel_Annals, Demailly_Kollar,Heier_ke_del_pezzo}). However, it is easier to handle \ke metrics instead of the Ricci flow with Nadel's method, because one can work with a $G$-invariant singular hermitian metric for $-K_X$ instead of $-2K_X$, resulting in a cohomology vanishing statement for a $G$-invariant multiplier ideal sheaf $\II$ of the form 
\begin{equation}\label{van_KE}
H^q (X,(K_X-K_X)\otimes \II)=H^q (X, \II)=0 \quad \forall q\geq 1.
\end{equation}
Note that \eqref{van_KE} yields more information on the zero-set of $\II$ than \eqref{van_RF}. In fact, the information in \eqref{van_KE} is strong enough to prove the existence of \ke metrics on all del Pezzo surfaces of degree $4,5,$ and $6$ (\cite{Heier_ke_del_pezzo}). In the present paper, we will show that the information in \eqref{van_RF} can be used to establish statement (ii) of Theorem \ref{PSS_Thm} for certain nongeneric del Pezzo surfaces with large automorphism group. In particular, we will prove\clearpage
\begin{maintheorem}\label{m_thm}
Let $X$ be one of the following del Pezzo surfaces.
\begin{enumerate}
\item $\PP^2$ blown up in four points in general position,
\item $\PP^2$ blown up in five points in general position with $\Aut(X)=\ZZ_2^4\rtimes\ZZ_4, \ZZ_2^4\rtimes(\ZZ_3\rtimes \ZZ_2),$ or $\ZZ_2^4\rtimes(\ZZ_5\rtimes \ZZ_2)$.
\end{enumerate}
Then statement (ii) of Theorem \ref{PSS_Thm} holds.
\end{maintheorem}
The method of proof also applies to certain del Pezzo surfaces of low degree, if their automorphism group is large enough, such as in the case of the Fermat cubic hypersurface in $\PP^3$ (see the remarks in Section \ref{low_degree}).
\begin{remark}
It is a result of Perelman that, given the existence of a \ke metric, the \kr flow will converge to it in the sense of Cheeger-Gromov. It should be noted that Theorem \ref{m_thm} does not assume the existence of a \ke metric. In fact, it proves the existence as an obvious corollary to the convergence statement (ii) of Theorem \ref{PSS_Thm}. Moreover, note again that the convergence in (ii) of Theorem \ref{PSS_Thm} is very strong, namely {\it exponentially fast in the $C^\infty$-norm}. 
\end{remark}

\begin{acknowledgement} The author would like to thank L.~Ein for helpful comments on an earlier version of this paper.
\end{acknowledgement}

\section{Classification and basic properties of del Pezzo surfaces}\label{delpezzoclass}
\begin{definition}
A {\it del Pezzo surface} is a two-dimensional compact complex manifold $X$ whose anti-canonical line bundle $-K_X$ is ample. We call the self-intersection number $(-K_X)^2=K_X^2$ the {\it degree} of $X$. We will denote the degree also by $d_X$.
\end{definition}
We now gather some important facts about del Pezzo surfaces, resulting in the standard classification (see \cite{Manin, Demazure, Hartshorne}).
\begin{facts}
For every del Pezzo surface $X$, the Picard group $\Pic X$ satisfies
\begin{equation*}
\rank \Pic X + d_X=10.
\end{equation*}
In particular, $d_X\leq 9$.\par
If $d_X=9$, then $X$ is isomorphic to $\PP^2$.\par
If $d_X=8$, then $X$ is isomorphic either to $\PP^1\times\PP^1$ or to $\tilde \PP^2$, ie $\PP^2$ blown up at one point.\par
If $7\geq d_X\geq 1$, then $X$ is isomorphic to $\PP^2$ blown up at $r=9-d_X$ points which have the following properties:
\begin{enumerate}
\item no three points lie on a line,
\item no six points lie on a conic,
\item no seven points lie on a cubic such that the eighth is a double point of the cubic.
\end{enumerate}
Any set of $r=9-d_X$ points satisfying the above three properties will be said to be in {\it general} position, and, conversely, the result of blowing up $1\leq r \leq 8$ points in general position in $\PP^2$ is a del Pezzo surface.. For $1\leq r\leq 4$ general points blown up, there is in each case a unique resulting del Pezzo surface. The reason is that for any two sets of points $P_1,\ldots,P_r$ and $Q_1,\ldots,Q_r$ ($r\leq 4$), with each set in general position, there is an element $A\in \Aut(\PP^2)=\PGL(3,\CC)$ with $A(P_i)=Q_i\ (1\leq i\leq r)$.
\end{facts}
For our understanding of del Pezzo surfaces, the following facts about the anti-canonical line bundle are also very important.
\begin{facts}
Let $1\leq r\leq 8$. Let $X$ be obtained by blowing up general points $P_i$, $i=1,\ldots,r$. Let $E_i$ denote the exceptional $(-1)$-curve that is the pre-image of $P_i$. Let $\pi:X\to\PP^2$ denote the blow up map. Then
\begin{equation*}
K_X=\pi^*K_{\PP^2}+\sum_{i=1}^{r}E_i.
\end{equation*}
This yields
\begin{equation*}
\dim H^0(X,-K_X)=10-r.
\end{equation*}\par
For $1\leq r\leq 6$, the complete linear system $|-K_X|$  gives an embedding into $\PP^{9-r}=\PP^{d_X}$. For $r=7$, it gives a double cover of $\PP^2$. The complete linear system $|-2K_X|$ gives an embedding into $\PP^6$. For $r=8$, $|-K_X|$ has a unique base point, $|-2K_X|$ gives a double cover of a singular quadric surface in $\PP^3$, and $|-3K_X|$ gives an embedding into $\PP^6$.\par

Finally, it turns out that, on every del Pezzo surface of degree at most $7$, the number of $(-1)$-curves exceeds $r$. The reason is that, when blowing up two points in $\PP^2$, the proper transform of the unique line through the two points becomes a $(-1)$-curve as well. When blowing up five points, the unique conic through the five points also becomes a $(-1)$-curve. It is easy to count these $(-1)$-curves: for $r=1,\ldots,8$, their numbers are $1,3,6,10,16,27,56,240$, respectively. Interestingly, for $r=1,\ldots,6$, under the map given by $|-K_X|$, all $(-1)$-curves become lines in projective space. Therefore, they are often referred to as {\it lines} on $X$. 
\end{facts}

\section{The case of 4 points blown up}
Let $X$ be a del Pezzo surface obtained from blowing up four points. We may, and do, assume these to be $P_1=[1,0,0]$, $P_2=[0,1,0]$, $P_3=[0,0,1]$, $P_4=[1,1,1]$. It is known that $\Aut(X)$ is the Weyl group of the root system of Dynkin type $D_5$, which is $S_5$ (see \cite{Koitabashi}, also \cite{Wiman}). However, we would like to understand $\Aut(X)$ more concretely.  \par
First all of, there is a subgroup $S_4$ of projectivities in $\PGL(3)=\Aut(\PP^2)$ that preserve the set $\{P_1,P_2,P_3,P_4\}$. These projectivities lift to $X$, and we can write $S_4\subset \Aut(X)$.\par
In addition, there exists for every $i=1,\ldots,4$ a quadratic Cremona transformation $\Cr_i$ that leaves $P_i$ fixed and has the three remaining points as indeterminacy locus. (Note that such a $\Cr_i$ is only defined up to the action of the $S_3\subset S_4$ consisting of automorphisms fixing $P_i$. For our purposes, it does not matter which $\Cr_i$ we choose.) All $\Cr_i$ extend to automorphisms of $X$. In light of this, we can write $\Aut(X)$ set-theoretically as a disjoint union
\begin{equation*}
\Aut(X)=S_4\uplus\left(\biguplus_{i=1}^4  \Cr_i \circ S_4\right).
\end{equation*}\par
In the following two Subsections, we will prove Theorem \ref{m_thm}, part (i), by means of Criterion \ref{simple_crit}. 
\subsection{Zero-dimensional multiplier ideal sheaves}
Let $\II\subseteq \OO_X$ be an $\Aut(X)$-invariant ideal sheaf satisfying Property (Van). In particular, for $E$ being the trivial line bundle,
\begin{equation}\label{van}
H^q(X,(-K_X)\otimes \II)=0 \quad \text{for all } q\geq 1.
\end{equation} 
Let $V=V(\II)$. Let $\dim V=0$, ie $V$ consists of a finite number of points. Consider the short exact sequence
\begin{equation*}
0\to \II(-K_X) \to \OO_X(-K_X)\to \OO_V(-K_X)\to 0.
\end{equation*}
Taking the corresponding long exact sequence yields
\begin{eqnarray*}
0&\to& H^0(X,\II(-K_X)) \to H^0(X,\OO_X(-K_X))\to H^0(V,\OO_V(-K_X))\to \\&\to& H^1(X,\II(-K_X)).
\end{eqnarray*}
From \eqref{van}, we have $H^1(X,\II(-K_X))=0$. Therefore, the map
\begin{equation*}
H^0(X,\OO_X(-K_X))\to H^0(V,\OO_V(-K_X))
\end{equation*}
is surjective. We saw in Section \ref{delpezzoclass} that $\dim H^0(X,\OO_X(-K_X))=10-4=6$. Therefore, $V$ consists of at most six points.
\begin{proposition}\label{r=4_no_0_dim}
There is no  $\Aut(X)$-invariant ideal sheaf $\II\subseteq \OO_X$ satisfying Property (Van) with $\dim V(\II)=0$. 
\end{proposition}
\begin{proof}
We assume that $\II$ exists and derive a contradiction. Since $V=V(\II)$ has at most six points, the contradiction arises from the claim that all orbits of $\Aut(X)$ have cardinality at least $8$. In the sequel, we prove this claim.\par
If the cardinality of an orbit were seven or less, then it would in fact be six or less, because the order of $\Aut(X)$ is even. The stabilizer subgroup of a point in an orbit of cardinality six or less would be of order $\frac{120}{6} =20$ or more. However, the only subgroups of $S_5$ of order $20$ or more are $S_5$, $A_5$, $S_4$, and the Frobenius group $\ZZ_5\rtimes\ZZ_4$, none of which has a faithful two-dimensional complex representation (on the tangent space to any point in the orbit). See eg \cite[\S 26]{Dornhoff} for more details.
\end{proof}
\subsection{One-dimensional multiplier ideal sheaves}
Let $\II\subseteq \OO_X$ be an $\Aut(X)$-invariant ideal sheaf satisfying Property (Van). Let $V=V(\II)$. Let $\dim_x V=1$ for all $x\in V$. We assume that the scheme defined by $\II\subseteq \OO_X$ has no embedded components.\par
It follows from Property (Van) and eg \cite[Theorem 2.1]{Nadel_Wuppertal} or \cite{PAGII} that
\begin{equation*}
R^i\pi_*(\II(-K_X))=0\ \text{for}\ i >0
\end{equation*}
and \begin{equation*}
H^i (\PP^2,\pi_*(\II(-K_X)))=0\ \text{for}\ i > 0.
\end{equation*}
In particular,
\begin{equation*}
H^2(\PP^2,\pi_*(\II(-K_X))=0.
\end{equation*}
By the standard projection formula (\cite[p. 281]{Fulton}) we have (for some $k_1,\ldots,k_4\in \NN$)
\begin{eqnarray*}
\pi_*(\II(-K_X))&=&\pi_*(\pi^*(-K_{\PP^2})\otimes\OO_{X}(-\sum_{i=1}^{4}E_i)\otimes\II)\\
&=&\OO_{\PP^2}(-K_{\PP^2})\otimes \mm_{P_1}^{k_1}\otimes\ldots\otimes\mm_{P_4}^{k_4}\otimes \fra,
\end{eqnarray*}
where $\fra$ is a locally principal ideal sheaf in $\OO_{\PP^2}$ whose zero-set consists of the images of the non-exceptional components of $V(\II)$.
Let $\fra\cong\OO_{\PP^2}(-d)$ with $d\geq 1$. Then by Serre duality:
\begin{eqnarray*}
0&=&H^2(\PP^2,\pi_*(\II(-K_X)))\\
&=& H^2(\PP^2,\OO_{\PP^2}(-K_{\PP^2})\otimes \mm_{P_1}^{k_1}\otimes\ldots\otimes\mm_{P_4}^{k_4}\otimes \fra)\\
&=& H^2(\PP^2,\OO_{\PP^2}(-K_{\PP^2})\otimes \fra)\\
&=& H^2(\PP^2,\OO_{\PP^2}(-K_{\PP^2})\otimes \OO_{\PP^2}(-d))\\
&=& H^0(\PP^2,\OO_{\PP^2}(2K_{\PP^2})\otimes \OO_{\PP^2}(d))\\
&=& H^0(\PP^2,\OO_{\PP^2}(d-6)).
\end{eqnarray*}
Note that this is the case if and only if $d\leq 5$. From this information, we will now extract information  on $V$ and show that in fact no such $V$ can exist. First we record two lemmas.
\begin{lemma}\label{r=4_eff_invar}
$\Aut(X)$ acts effectively on any $\Aut(X)$-invariant irreducible curve.
\end{lemma}
\begin{proof}
For any given element of $\Aut(X)$, it is easy to list the irreducible curves that are left point-wise fixed by the given element (if any exist). However, none of these curves are $\Aut(X)$-invariant. 
\end{proof}
\begin{lemma}[{\cite[Excercise 20.18]{Harris}}]\label{num_sings}
The number of singular points of an irreducible plane curve of degree $d$ is no more than its arithmetic genus, which is $g_a = \frac 1 2 (d-1)(d-2)$.
\end{lemma}
\begin{proposition}\label{r=4_no_1_dim}
There is no  $\Aut(X)$-invariant ideal sheaf $\II\subseteq \OO_X$ satisfying Property (Van) with $\dim_x V(\II)=1$ for all $x\in V(\II)$ and such that the scheme defined by $\II\subseteq \OO_X$ has no embedded components. 
\end{proposition}
\begin{proof}
We assume that $\II$ exists and derive a contradiction. We can assume w.l.o.g.\ that no $(-1)$-curve is contained in $V=V(\II)$, because as a consequence of the $\Aut(X)$-invariance, all ten $(-1)$-curves would otherwise be contained in $V$, and $\pi(V)$ would have at least six irreducible one-dimensional components, in violation of $H^0(\PP^2,\OO_{\PP^2}(d-6))=0$.\par
Since all elements $f$ of $\Aut(X)$ are induced by an automorphism or birational map of $\PP^2$ (again denoted by $f$), we note that for an irreducible curve $D$ not contained in the exceptional set of $\pi$, $\pi(f(D))=f(\pi(D))$, where $f(\pi(D))$ is defined to be the closure of $f(\pi(D)\cap\Dom(f))$.\par
First, let us assume that $V$ is irreducible.\par
If $d=1,2$, then $V=\PP^1$. This is a contradiction, because $S_5\subset \Aut(X)$ acts effectively on $V$ by Lemma \ref{r=4_eff_invar}, but $S_5\not\subseteq \PGL(2,\CC)$ (see \cite[Chapter 2]{Grove_Benson}).\par
If $d=3$, then $g_a=1$. Since $\Aut(X)$ acts on $V$ with orbits of length at least $8$, $V$ is smooth by Lemma \ref{num_sings}. We saw in Lemma \ref{r=4_eff_invar} that $S_5$ acts effectively on all $S_5$-invariant curves. However, it is not a subgroup of the automorphism group of any elliptic curve, because it is well known (see \cite[p.~64]{Miranda}) that the possible automorphism groups of an elliptic curve are 
\begin{equation}\label{aut_torus}
\ZZ_2\ltimes\CC/\Gamma,\ZZ_4\ltimes\CC/\Gamma,\ZZ_6\ltimes\CC/\Gamma,
\end{equation}
none of which can contain $S_5$. We have obtained a contradiction.\par
If $d=4$, then $g_a=3$. Again, $V$ is smooth. As an abstract group, the $2$-Sylow subgroup of $S_5$ is $\ZZ_4\rtimes\ZZ_2$, which is given concretely eg by the subgroup
\begin{equation*}
\{\id,(12),(34),(12)(34),(14)(23),(13)(24),(1324),(1423)\}\subset S_5.
\end{equation*}
The group $\ZZ_4\rtimes\ZZ_2$ has cyclic $2$-deficiency equal to $1$, as defined in \cite[p.197]{Kulkarni}, and by \cite[Theorem 2.3]{Kulkarni} $S_5$ does not act on $V$. A contradiction.\par
If $d=5$, then $g_a=6$. Again, $V$ is smooth, and another application of \cite[Theorem 2.3]{Kulkarni} yields again a contradiction. \par
Next, we treat the case where the number of irreducible components of $V$ is $2$.\par
For $d=2$, the two components of $\pi(V)$ are lines. If the two lines intersected outside of $\{P_1,P_2,P_3,P_4\}$, then the components of $V$ would intersect in precisely one point as well. However, this is impossible, because the minimum orbit length of the action of $\Aut(X)$ is $8$. Therefore, the two lines must each go through the same $P_i$. It is now obvious that they are not invariant under $S_4\subset S_5$.\par
For $d=3$, the components of $\pi(V)$ are a line and a conic. They must meet in one (with multiplicity two) or two of the $\{P_1,P_2,P_3,P_4\}$. We again obtain a contradiction.\par
For $d=4$, there are two cases. First, the components of $V$ could be a rational curve and an elliptic curve. They each have to be invariant under $S_5$, in which case we obtained a contradiction earlier. Second, the components of $\pi(V)$ could be two conics. Since the length of any orbit of the action of $S_5$ is at least $8$, these conics must intersect precisely in the points $P_1=[1,0,0]$, $P_2=[0,1,0]$, $P_3=[0,0,1]$, $P_4=[1,1,1]$, so that they get separated under the blowing up. Note that any conic through these four points is of the form
\begin{equation*}
a_0 X_1X_2+ a_1 X_0X_2+(-a_0-a_1)X_0X_1=0,
\end{equation*}
and any such conic is mapped to a line under $\Cr_4$, a contradiction.\par
For $d=5$, the components of $V$ are a rational curve and an elliptic curve or a rational curve and a curve of genus $3$. From the rational curve, which is preserved, we get a contradiction.\par
Now we treat the case of three irreducible components.\par
For $d=3$, the components of $\pi(V)$ are lines, meeting in at most three of the points $P_1,P_2,P_3,P_4$. Thus they cannot be invariant under the $S_4$ action.\par
For $d=4$, the three components of $\pi(V)$ are two lines and one conic. All points of intersection must be contained in $\{P_1,P_2,P_3,P_4\}$, but the multiplicity will not be the same at all $P_i$, $i=1,\ldots,4$. This is again impossible by the $S_4$-symmetry.\par
For $d=5$, $\pi(V)$ consists of either two lines and an elliptic curve or one line and two conics. In either case, we can argue as we did above and obtain a contradiction.\par
Now for the case of four irreducible components.\par
If $d=4$, the components of $\pi(V)$ are lines. Their at most six points of intersection must be $\{P_1,P_2,P_3,P_4\}$ (counted with multiplicity). However, the multiplicity cannot be distributed symmetrically, and we get a contradiction from the action of $S_4$.\par
If $d=5$, then the components of $\pi(V)$ are three lines and a conic. The eight points of intersection (counted with multiplicity) must either be contained entirely in $\{P_1,P_2,P_3,P_4\}$ or disjoint from it. Containment is impossible because of the $S_4$-symmetry. Now, again by the $S_4$-symmetry, all three lines must be disjoint from $\{P_1,P_2,P_3,P_4\}$. As a consequence, the equations of the three lines are of the form
\begin{equation*}
a_0X_0+a_1X_1+a_2X_2=0
\end{equation*}
with $a_0,a_1,a_2\not =0$. Under $\Cr_4$ such lines are mapped to conics of the form
\begin{equation*}
a_0X_1X_2+a_1X_0X_2+a_2X_0X_1=0,
\end{equation*}
which is a contradiction.\par
Finally, we treat the case of five components.\par
The only possibility is $d=5$ and $\pi(V)$ consisting of five lines. Since lines that are disjoint from $\{P_1,P_2,P_3,P_4\}$ are mapped to conics by $\Cr_4$, all five lines must have nonempty intersection with $\{P_1,P_2,P_3,P_4\}$. Clearly, one point must be contained in two lines. By the $S_4$-symmetry, all points must be contained in two lines. Arguing again with the minimum orbit length of $8$, we see that it is impossible for six or less points of intersection to be outside $\{P_1,P_2,P_3,P_4\}$, so all ten are in this set. However, it is impossible to distribute the ten points onto the four points with equal multiplicity. Therefore, the $S_4$-action gives our final contradiction.
\end{proof}
Propositions \ref{r=4_no_0_dim} and \ref{r=4_no_1_dim} prove Theorem \ref{m_thm}, part (i). Notice that we can assume w.l.o.g. that $V$ has the same dimension at all its points. Also, the assumption of nonexistence of embedded components in Proposition \ref{r=4_no_1_dim} is not a problem, since such embedded components can be ruled out by an $H^0$ computation similar to the one leading to the proof of Proposition \ref{r=4_no_0_dim}.

\section{Five points blown up}\label{r=5}
Let $X$ be a del Pezzo surface obtained by blowing up five points. We can find an automorphism of $\PP^2$ that takes the five points to $P_1=[1,0,0]$, $P_2=[0,1,0]$, $P_3=[0,0,1]$, $P_4=[1,1,1]$, $P_5=[a,b,c]$, with $(a,b,c) \in(\CC^*\times\CC^*\times\CC^*)\backslash\{(1,1,1)\}$. (The reason for $a,b,c \not= 0$ is that no three of these point lie on a line.)\par
The structure of $\Aut(X)$ is described eg in \cite{Hosoh_quartic} (see also \cite{Wiman}). It turns out that it is always of the form
\begin{equation*}
\Aut(X)=\ZZ_2^4\rtimes G_{P_5},
\end{equation*}
where $G_{P_5}$ is a subgroup of $S_5$ depending on the point $P_5$. The possibilities for $G_{P_5}$ are
\begin{enumerate}
\item $\{\id\},$
\item $\ZZ_2,$
\item $\ZZ_4,$
\item $\ZZ_3\rtimes \ZZ_2,$
\item $\ZZ_5\rtimes \ZZ_2.$
\end{enumerate}
The elements of $G_{P_5}$ are lifts of those elements of $\PGL(3,\CC)$ that map the set $\{P_1,P_2,P_3,P_4,P_5\}$ to itself. For a generic point $P_5$, we have $G_{P_5}=\{\id\}$. More precisely:
\begin{proposition}\label{GP5_nontrival}
One has $G_{P_5}\not=\{id\}$ if and only if $P_5=[1,\xi,1+\xi]$ with $\xi \in \CC\backslash\{0,1,-1\}$. Moreover, $G_{P_5}=\ZZ_2$ holds precisely when $\xi^2+1\not = 0$ and $\xi^2\pm\xi\pm 1\not =0$.
\end{proposition}\par
When $G_{P_5}\not=\{id\}$, \cite[Proposition 8.1.11]{Blanc_thesis_lanl} gives explicitly the elements of $G_{P_5}\subset \PGL(3,\CC)$ and their action on the set $\{P_1,P_2,P_3,P_4,P_5\}$. Based on Proposition \ref{GP5_nontrival}, it is clear that  there is one element of $\PGL(3,\CC)$ that is contained in $G_{P_5}$ whenever it is not the trivial group, namely 
\begin{equation*}
[X_0,X_1,X_2]\mapsto [X_2-X_1,X_2-X_0,X_2].
\end{equation*}\par
Next, let us have a closer look at the elements of $\ZZ_2^4\subseteq \Aut(X)$. The following two birational involutions of $\PP^2$ lift to elements of $\Aut(X)$.\par
We define $\Cr_{45}$  to be
\begin{equation*}
\Cr_{45}([X_0,X_1,X_2])=[aX_1X_2,bX_0X_2,cX_0X_1].
\end{equation*}
This is a quadratic Cremona transformation that exchanges $P_4$ and $P_5$ and has $\{P_1,P_2,P_3\}$ as indeterminacy locus. We abuse notation and refer both to the birational involution and the corresponding element of $\Aut(X)$ with the symbol $\Cr_{45}$. \par
Moreover, we let $\sigma_1$ be the following cubic birational involution of $\PP^2$.
\begin{eqnarray*}
\sigma_1([X_0,X_1,X_2])&=&[-aX_1X_2((c-b)X_0+(a-c)X_1+(b-a)X_2),\\
&&X_1(a(c-b)X_1X_2+b(a-c)X_0X_2+c(b-a)X_0X_1),\\
&&X_2(a(c-b)X_1X_2+b(a-c)X_0X_2+c(b-a)X_0X_1)].
\end{eqnarray*}
We shall refer both to the birational involution and the corresponding element of $\Aut(X)$ with the symbol $\sigma_1$. It is easy to see that the strict transform of the cubic curve $C$ given by
\begin{equation*}
b(a-c)X_0^2X_2+c(b-a)X_0^2X_1+a(a-c)X_1^2X_2+a(b-a)X_1X_2^2+2a(c-b)X_0X_1X_2=0
\end{equation*}
is precisely the set of points fixed point-wise by the lift of above $\sigma_1$. In addition, the following Lemma tells us that $C$ is invariant under every element of $\ZZ_2^4$.
\begin{lemma}
Let $A$ be an abelian group acting on a set $M$. For $g\in A$, let 
\begin{equation*}
M^g=\{ x\in M: gx=x\}.
\end{equation*}
Then $M^g$ is invariant under every element of $A$.
\end{lemma}
\begin{proof}
For $x\in M^g$, we have for any $h\in A$:
\begin{equation*}
g(hx)=h(gx)=hx, 
\end{equation*}
ie $hx\in M^g$ also.
\end{proof}
Closer inspection (\cite{Wiman}, \cite{Koitabashi}, \cite{Hosoh_quartic}) shows that $\ZZ_2^4\subseteq\Aut(X)$ contains the lifts of ten quadratic Cremona involutions $\Cr_{ij}\ (1\leq i<j\leq 5)$ that exchange $P_i,P_j$ and have the remaining three points as indeterminacy locus.  Again, we abuse notation and denote the maps before and after the lift by the same symbols $\Cr_{ij}$.\par

Moreover, $\ZZ_2^4\subseteq\Aut(X)$ contains the lifts of five cubic involutions $\sigma_i\ (1\leq i\leq 5)$. By comparing the respective action on the set of the sixteen $(-1)$-curves (which determines any automorphism uniquely, see \cite{Koitabashi} or \cite{Hosoh_quartic}), it is easily verified that 
\begin{eqnarray*}
\Cr_{ij}\circ\Cr_{kl}&=&\Cr_{kl} \circ\Cr_{ij},\\
\Cr_{ij}\circ\Cr_{jk}&=&\Cr_{ik}.
\end{eqnarray*}
Moreover, for $i,j,k,l,m$ all distinct, we have
\begin{equation*}
\sigma_m=\Cr_{ij}\circ\Cr_{kl}.
\end{equation*}
In particular, we have
\begin{equation}\label{group_law}
\sigma_j=\Cr_{1j}\circ\sigma_1 \quad (2\leq j\leq 5).
\end{equation}
\par
The $\Cr_{ij}$ have precisely four fixed points each (both before and after the lifting). Thus, the only $\ZZ_2^4$-invariant curves on which $\ZZ_2^4$ does not act effectively are the lifts of the curves $C$ pertaining to the cubic involutions. However, for $a,b,c\not =0$, the curves $C$ are smooth elliptic curves (easy exercise). Therefore, the group $\ZZ_2^4\subset \Aut(X)$ acts effectively on any $\Aut(X)$-invariant irreducible curve that is not an elliptic curve.\par
In the following two Subsections, we will prove Theorem \ref{m_thm}, part (ii), by means of Criterion \ref{simple_crit}.
\subsection{Zero-dimensional multiplier ideal sheaves}
Let $\II\subseteq \OO_X$ be an $\Aut(X)$-invariant ideal sheaf satisfying Property (Van). Let $V=V(\II)$. Let $\dim V=0$. 
We saw in the previous section that the map
\begin{equation*}
H^0(X,\OO_X(-K_X))\to H^0(V,\OO_V(-K_X))
\end{equation*}
is surjective. Since $\dim H^0(X,\OO_X(-K_X))=10-5=5$, $V$ consists of at most five points.
\begin{proposition}\label{r=5_no_0_dim}
If $G_{P_5}=\ZZ_4, \ZZ_3\rtimes \ZZ_2,$ or $\ZZ_5\rtimes \ZZ_2$, then there is no $\Aut(X)$-invariant ideal sheaf $\II\subseteq \OO_X$ satisfying Property (Van) with $\dim V(\II)=0$. 
\end{proposition}
\begin{proof}
We assume that $\II$ exists and derive a contradiction. Since $V=V(\II)$ has at most five points, the statement follows from the claim that all orbits of $\Aut(X)$ have cardinality at least $8$. In the sequel, we prove this claim.\par
We first consider the action of $\ZZ_2^4\subseteq \Aut(X)$ on $X$. There are no orbits of length less than $4$, which can be shown as follows: if we assume that there is such an orbit, then the stabilizer subgroup $(\ZZ_2^4)_P$ of a point $P$  in the orbit would be a subgroup of order $8$ or $16$ in $\ZZ_2^4$. The only such groups are $\ZZ_2^3$ and $\ZZ_2^4$, but by Schur's Lemma these groups do not permit faithful two-dimensional complex representations (on the tangent space to any point in the orbit), contradiction.\par
Now let $P\in X$ be such that the cardinality of the orbit of $P$ under $\ZZ_2^4$ is $\# \ZZ_2^4 P = 4$. The list in \cite[Proposition 8.1.11]{Blanc_thesis_lanl} gives the action of $G_{P_5}$ on the elements of the set $\{P_1,P_2,P_3,P_4,P_5\}$. This data determines the action on the sixteen $(-1)$-curves and therefore determines the automorphism uniquely.\par
When $G_{P_5}=\ZZ_4, \ZZ_3\rtimes \ZZ_2,$ or $\ZZ_5\rtimes \ZZ_2$, the claim is implied by $\# (\ZZ_2^4\rtimes G_{P_5})_P \leq 8$. A case by case analysis shows that if $\# (\ZZ_2^4\rtimes G_{P_5})_P > 8$, then 
\begin{equation*}
\ZZ_2^3\subseteq (\ZZ_2^4\rtimes G_{P_5})_P \cap \ZZ_2^4,
\end{equation*}
which yields a contradiction by Schur's Lemma. Note that the strict inequality $\# (\ZZ_2^4\rtimes G_{P_5})_P > 8$ is necessary, because otherwise $(\ZZ_2^4\rtimes G_{P_5})_P=\ZZ_2^2\rtimes \ZZ_2$ might (and does) occur.
\end{proof}

\subsection{One-dimensional multiplier ideal sheaves}
In this subsection, we prove
\begin{proposition}\label{r=5_no_1_dim}
If $G_{P_5}=\ZZ_4, \ZZ_3\rtimes \ZZ_2,$ or $\ZZ_5\rtimes \ZZ_2$, then there is no  $\Aut(X)$-invariant ideal sheaf $\II\subseteq \OO_X$ satisfying Property (Van) with $\dim_x V(\II)=1$ for all $x\in V(\II)$ and such that the scheme defined by $\II\subseteq \OO_X$ has no embedded components. 
\end{proposition}
\begin{proof}
We assume that $\II$ exists and derive a contradiction. We can again assume w.l.o.g. that no $(-1)$-curve is contained in $V=V(\II)$.\par
First, let us assume that $V$ is irreducible. \par
If $d=1,2$, then $V=\PP^1$. This is a contradiction, because $\ZZ_2^4\subset \Aut(X)$ acts effectively on $V$, but $\ZZ_2^4\not\subseteq \PGL(2,\CC)$ (see \cite[Chapter 2]{Grove_Benson}).\par
If $d=3$, then $g_a=1$. Since $\Aut(X)$ acts on $V$ with orbits of length at least $8$, $V$ is smooth by Lemma \ref{num_sings}. We know that $\ZZ_2^3\rtimes G_{P_5}$ acts effectively on $V$. However, it is not a subgroup of any of the groups listed in \eqref{aut_torus}. We have obtained a contradiction.\par
If $d=4$, then $g_a=3$. Again, $V$ is smooth. We can obtain a contradiction both by Maclachlan's theorem on abelian groups of automorphisms of Riemann surfaces \cite[Theorem 4]{Maclachlan} and by the previously used \cite[Theorem 2.3]{Kulkarni}, noting that the cyclic $2$-deficiency of $\ZZ_2^4$ is equal to $3$.\par
If $d=5$, then $g_a=6$. Again, $V$ is smooth. \cite[Theorem 2.3]{Kulkarni} applies also in this case, and we obtain a contradiction. \par
Next, we treat the case where the number of irreducible components is $2$.\par
For $d=2,3$, the two components are smooth rational curves. If each component (call them $V_1,V_2$) is invariant under $\ZZ_2^4\subset\Aut(X)$, we have a contradiction as before.\par
Let us assume that there is $g\in\ZZ_2^4\subset \Aut(X)$ such that $g(V_1)=V_2$.  Let
\begin{equation*}
G_1=\{g\in \ZZ_2^4: g(V_1)=V_1\}
\end{equation*}
Then 
\begin{equation*}
\ZZ_2^4 =G_1\uplus gG_1.
\end{equation*}
Therefore, the index of $G_1$ in $\ZZ_2^4$ is two, and consequently $G_1=\ZZ_2^3$, which gives a contradiction, since $\ZZ_2^3\not\subseteq \PGL(2,\CC)$.\par
When $d=4$, $V$ may be the disjoint union of a smooth elliptic curve and a rational curve or the disjoint union of two smooth rational curves. In either case, this is impossible.\par
When $d=5$, $V$ may be the disjoint union of a smooth elliptic curve and a smooth rational curve or the disjoint union of a smooth curve of genus $3$ and a smooth rational curve. In any case, this is impossible.\par
Next, let us assume that there are three irreducible components $V_1,V_2,V_3$.\par
When $d=3,4$, $V$ is the union of three smooth rational curves. If the action of $\ZZ_2^4\subset \Aut(X)$ on the three components leaves one component invariant, then $\ZZ_2^4\subset \PGL(2,\CC)$, which is a contradiction.\par
If no component is invariant under $\ZZ_2^4\subset \Aut(X)$, then there are $g_2,g_3\in \ZZ_2^4$ such that $g_2(V_1)=V_2, g_3(V_1)=V_3$. However, a brief computation reveals $(g_2\cdot g_3)^2(V_1)\in V_2$, which means that the order of $g_2\cdot g_3$ cannot be $2$. However, all nontrivial elements of $\ZZ_2^4$ have order two, a contradiction.\par
When $d=5$, one of the components might be a smooth elliptic curve, but in any case one obtains a contradiction as before.\par

Finally, in the cases of $4$ and $5$ irreducible components, all components are rational curves. If there are two components whose union is $\ZZ_2^4$ invariant, then we obtain a contradiction as in the case of two irreducible componenets. If not, then there exist $g_2,g_3\in \ZZ_2^4$ and components $V_1,V_2,V_3$ such that $g_2(V_1)=V_2, g_3(V_1)=V_3$. Again, we obtain a contradiction.
\end{proof}

Propositions \ref{r=5_no_0_dim} and \ref{r=5_no_1_dim} prove Theorem \ref{m_thm}, part (ii).

\section{Comments on the cases of $r\not =4,5$}\label{low_degree}

\subsection{One or two points blown up}
For $\PP^2$ blown up at one or two points, one can show that the so-called Calabi-Futaki invariant does not vanish (see eg \cite[Examples 3.10, 3.11]{Tian_book} for details). The nonvanishing of this invariant is an obstruction to the existence of a \ke metric. If $g_{i\bar{j}}(t)$ ($t\to \infty)$ were to converge, it would necessarily converge against a \ke metric. So the statement of Theorem \ref{PSS_Thm} (ii) cannot hold on 
$\PP^2$ blown up at one or two points.
\subsection{Three points blown up}
In this paper, we do not make a statement about this case, because both zero- and one-dimensional zero-sets of multiplier ideal sheaves cannot be ruled out using Property (Van).\par
In the zero-dimensional case, the problem is that the surjectivity of the map
\begin{equation*}
H^0(X,\OO_X(-K_X))\to H^0(V,\OO_V(-K_X))
\end{equation*}
limits the cardinality of the zero-set to $10-3=7$. However, the six points of intersection in the ``hexagon'' formed by the six $(-1)$-curves are clearly $\Aut(X)$-invariant, and we are unable to rule out that they are the $\Aut(X)$-invariant zero-set of a multiplier ideal sheaf satisfying \eqref{van}.\par
Similarly, the union of the six $(-1)$-curves is clearly $\Aut(X)$-invariant, and it is not possible to rule out that it forms the zero-set of a multiplier ideal sheaf based on Property (Van). Clearly, this case merits further investigation.

\subsection{Six or more points blown up}
For a generic del Pezzo surface $X$ of degree three, $\Aut(X)$ is unfortunately the trivial group (see \cite{Koitabashi}). There are of course non-generic $X$ which have extra automorphisms, and a nice list of these can be found in \cite[Table 10.3]{Dolgachev}. Recall that del Pezzo surfaces of degree three are precisely the smooth cubic hypersurfaces of $\PP^3$, and perhaps the most important example is the Fermat cubic surface in $\PP^3$ given by
\begin{equation*}
Z_0^3+Z_1^3+Z_2^3+Z_3^3=0.
\end{equation*}
Its automorphism group is $\ZZ_3^3\rtimes S_4$, acting in the obvious way. The cardinality of this group is $27\cdot 24=648$. In this case, the same analysis as in Section \ref{r=5} does yield statement (ii) of Theorem \ref{PSS_Thm}. In fact, one can even do without the results concerning the nature of the automorphism group for Riemann surfaces of genus no more than 6, because according to the well-known Hurwitz bound, such an automorphism group has cardinality at most $84(g-1)$, which is less than $648$. We leave the details to the reader.\par
At the other end of the spectrum, on a del Pezzo surface $X$ of degree one, the unique base point of the linear system $|-K_X|$ is fixed by all automorphisms. Therefore, unfortunately, $\Aut(X)$ acts with a fixed point regardless of the nature of $X$, and we are unable to handle to this case.\par
On a del Pezzo surface $X$ of degree two, the linear system $|-K_X|$ gives a two-sheeted cover of $\PP^2$ branched along a smooth curve $C$ of degree $4$ in $\PP^2$. This cover defines an involutive automorphism of $X$ called the Geiser involution. On a generic $X$, this is the only nontrivial automorphism, ie $\Aut(X)=\ZZ_2$. However, certain non-generic $X$ do have extra automorphisms. A list of these $X$ and their automorphism groups, together with a lucid exposition of the topic, can be found in \cite[Table 10.4]{Dolgachev}. We do not go into any details regarding this case.


\begin{thebibliography}{Wim96}

\bibitem[Bla06]{Blanc_thesis_lanl}
J.~Blanc.
\newblock Finite abelian subgroups of the {C}remona group of the plane.
\newblock {\em {\rm arXiv:math.AG/0610368}}, 2006.

\bibitem[Cao85]{Cao_Inv_Math_1985}
H.-D. Cao.
\newblock Deformation of {K}\"ahler metrics to {K}\"ahler-{E}instein metrics on
  compact {K}\"ahler manifolds.
\newblock {\em Invent. Math.}, 81(2):359--372, 1985.

\bibitem[CCZ03]{Cao_Chen_Zhu}
H.-D. Cao, B.-L. Chen, and X.-P. Zhu.
\newblock Ricci flow on compact {K}\"ahler manifolds of positive bisectional
  curvature.
\newblock {\em C. R. Math. Acad. Sci. Paris}, 337(12):781--784, 2003.

\bibitem[CT06]{Chen_Tian}
X.~X. Chen and G.~Tian.
\newblock Ricci flow on {K}\"ahler-{E}instein manifolds.
\newblock {\em Duke Math. J.}, 131(1):17--73, 2006.

\bibitem[Dem80]{Demazure}
M.~Demazure.
\newblock Surfaces de {D}el {P}ezzo, {I}-{V}.
\newblock {\em {L}ecture {N}otes in {M}athematics ({S}pringer-{V}erlag)},
  777:21--61, 1980.

\bibitem[DK01]{Demailly_Kollar}
J.-P. Demailly and J.~Koll{\'a}r.
\newblock Semi-continuity of complex singularity exponents and
  {K}\"ahler-{E}instein metrics on {F}ano orbifolds.
\newblock {\em Ann. Sci. \'Ecole Norm. Sup. (4)}, 34(4):525--556, 2001.

\bibitem[Dol07]{Dolgachev}
I.~V. Dolgachev.
\newblock Topics in classical algebraic geometry.
\newblock {\em {L}ecture notes, available at {\tt
  http://www.math.lsa.umich.edu/{$\tilde{\ }$}idolga/lecturenotes.html}}, 2007.

\bibitem[Dor71]{Dornhoff}
L.~Dornhoff.
\newblock {\em Group representation theory. {P}art {A}: {O}rdinary
  representation theory}.
\newblock Marcel Dekker Inc., New York, 1971.
\newblock Pure and Applied Mathematics, 7.

\bibitem[FH91]{Fulton_Harris}
W.~Fulton and J.~Harris.
\newblock {\em Representation theory}, volume 129 of {\em Graduate Texts in
  Mathematics}.
\newblock Springer-Verlag, New York, 1991.
\newblock A first course, Readings in Mathematics.

\bibitem[Ful98]{Fulton}
W.~Fulton.
\newblock {\em Intersection theory}, volume~2 of {\em Ergebnisse der Mathematik
  und ihrer Grenzgebiete. 3. Folge. A Series of Modern Surveys in Mathematics}.
\newblock Springer-Verlag, Berlin, second edition, 1998.

\bibitem[GB85]{Grove_Benson}
L.~C. Grove and C.~T. Benson.
\newblock {\em Finite reflection groups}, volume~99 of {\em Graduate Texts in
  Mathematics}.
\newblock Springer-Verlag, New York, second edition, 1985.

\bibitem[Har77]{Hartshorne}
R.~Hartshorne.
\newblock {\em Algebraic geometry}.
\newblock Springer-Verlag, New York, 1977.
\newblock Graduate Texts in Mathematics, No. 52.

\bibitem[Har92]{Harris}
J.~Harris.
\newblock {\em Algebraic geometry}, volume 133 of {\em Graduate Texts in
  Mathematics}.
\newblock Springer-Verlag, New York, 1992.

\bibitem[Hei07]{Heier_ke_del_pezzo}
G.~Heier.
\newblock Existence of {K}\"ahler-{E}instein metrics and multiplier ideal
  sheaves on del {P}ezzo surfaces.
\newblock {\em {\rm arXiv:math/0710.5724}}, 2007.

\bibitem[Hos96]{Hosoh_quartic}
T.~Hosoh.
\newblock Automorphism groups of quartic del {P}ezzo surfaces.
\newblock {\em J. Algebra}, 185(2):374--389, 1996.

\bibitem[Hos97]{Hosoh_cubic}
T.~Hosoh.
\newblock Automorphism groups of cubic surfaces.
\newblock {\em J. Algebra}, 192(2):651--677, 1997.

\bibitem[Koh79]{Kohn_Acta}
J.~J. Kohn.
\newblock Subellipticity of the {$\bar \partial $}-{N}eumann problem on
  pseudo-convex domains: sufficient conditions.
\newblock {\em Acta Math.}, 142(1-2):79--122, 1979.

\bibitem[Koi88]{Koitabashi}
M.~Koitabashi.
\newblock Automorphism groups of generic rational surfaces.
\newblock {\em J. Algebra}, 116(1):130--142, 1988.

\bibitem[Kul87]{Kulkarni}
R.~Kulkarni.
\newblock Symmetries of surfaces.
\newblock {\em Topology}, 26(2):195--203, 1987.

\bibitem[Laz04]{PAGII}
R.~Lazarsfeld.
\newblock {\em Positivity in algebraic geometry. {II}}, volume~49 of {\em
  Ergebnisse der Mathematik und ihrer Grenzgebiete. 3. Folge. A Series of
  Modern Surveys in Mathematics}.
\newblock Springer-Verlag, Berlin, 2004.
\newblock Positivity for vector bundles, and multiplier ideals.

\bibitem[Mac65]{Maclachlan}
C.~Maclachlan.
\newblock Abelian groups of automorphisms of compact {R}iemann surfaces.
\newblock {\em Proc. London Math. Soc. (3)}, 15:699--712, 1965.

\bibitem[Mac90]{Macbeath}
A.~M. Macbeath.
\newblock Automorphisms of {R}iemann surfaces.
\newblock In {\em Combinatorial group theory (College Park, MD, 1988)}, volume
  109 of {\em Contemp. Math.}, pages 107--112. Amer. Math. Soc., Providence,
  RI, 1990.

\bibitem[MH74]{Manin}
Yu.~I. Manin and M.~Hazewinkel.
\newblock {\em Cubic forms: algebra, geometry, arithmetic}.
\newblock North-Holland Publishing Co., Amsterdam, 1974.
\newblock Translated from Russian by M. Hazewinkel, North-Holland Mathematical
  Library, Vol. 4.

\bibitem[Mir95]{Miranda}
R.~Miranda.
\newblock {\em Algebraic curves and {R}iemann surfaces}, volume~5 of {\em
  Graduate Studies in Mathematics}.
\newblock American Mathematical Society, Providence, RI, 1995.

\bibitem[Nad89]{Nadel_PNAS}
A.~M. Nadel.
\newblock Multiplier ideal sheaves and existence of {K}\"ahler-{E}instein
  metrics of positive scalar curvature.
\newblock {\em Proc. Nat. Acad. Sci. U.S.A.}, 86(19):7299--7300, 1989.

\bibitem[Nad90]{Nadel_Annals}
A.~M. Nadel.
\newblock Multiplier ideal sheaves and {K}\"ahler-{E}instein metrics of
  positive scalar curvature.
\newblock {\em Ann. of Math. (2)}, 132(3):549--596, 1990.

\bibitem[Nad91]{Nadel_Wuppertal}
A.~M. Nadel.
\newblock The behavior of multiplier ideal sheaves under morphisms.
\newblock In {\em Complex analysis (Wuppertal, 1991)}, Aspects Math., E17,
  pages 205--222. Vieweg, Braunschweig, 1991.

\bibitem[PS07]{PS_CDM}
D.~H. Phong and J.~Sturm.
\newblock Lectures on stability and constant scalar curvature.
\newblock {\em {L}ecture notes for the conference {\rm Current Developments in
  Mathematics, Harvard University}}, 2007.

\bibitem[PSS06]{PSS}
D.~H. Phong, N.~Sesum, and J.~Sturm.
\newblock Multiplier ideal sheaves and the {K}\"ahler-{R}icci flow.
\newblock {\em {\rm arXiv:math/0611794v2}}, 2006.

\bibitem[Tia90]{Tian_Inv_Math}
G.~Tian.
\newblock On {C}alabi's conjecture for complex surfaces with positive first
  {C}hern class.
\newblock {\em Invent. Math.}, 101(1):101--172, 1990.

\bibitem[Tia00]{Tian_book}
G.~Tian.
\newblock {\em Canonical metrics in {K}\"ahler geometry}.
\newblock Lectures in Mathematics ETH Z\"urich. Birkh\"auser Verlag, Basel,
  2000.
\newblock Notes taken by Meike Akveld.

\bibitem[TZ07]{Tian_Zhu}
G.~Tian and X.~Zhu.
\newblock Convergence of {K}\"ahler-{R}icci flow.
\newblock {\em J. Amer. Math. Soc.}, 20(3):675--699 (electronic), 2007.

\bibitem[Wim96]{Wiman}
A.~Wiman.
\newblock Zur {T}heorie der endlichen {G}ruppen von birationalen
  {T}ransformationen in der {E}bene.
\newblock {\em Math. Ann.}, 48(1-2):195--240, 1896.

\end{thebibliography}
\end{document}